\documentclass{article}

\usepackage{multirow}
\usepackage{epsf}
\usepackage{latexsym}
\usepackage{epsfig}
\usepackage{natbib}
\usepackage{tabularx}
\usepackage{array}
\usepackage{amsmath, amssymb, amsthm}
\usepackage[ruled, vlined]{algorithm2e}
\usepackage{multicol}
\usepackage{comment}
\usepackage{lscape}
\usepackage{natbib}

\newtheorem{thm}{Theorem}
\newtheorem*{pf}{Proof}

\title{Confidence-based Optimization for the Newsvendor Problem}
\author{Roberto Rossi\\ 
{\small University of Edinburgh Business School, UK,}\\{\small roberto.rossi@ed.ac.uk}\\
Steven Prestwich\\
{\small Dept. of Computer Science, University College Cork, Ireland,}\\{\small s.prestwich@4c.ucc.ie}\\
S. Armagan Tarim\\
{\small Dept. of Management, Hacettepe University, Turkey,}\\{\small armagan.tarim@hacettepe.edu.tr}\\
Brahim Hnich\\
{\small Dept. of Computer Engineering, Izmir University of Economics, Turkey,}\\{\small brahim.hnich@ieu.edu.tr}}

\date{}

 \DeclareMathOperator*{\argmin}{arg\,min}
  \DeclareMathOperator*{\argmax}{arg\,max}

\begin{document}
\parskip 1pc 
\setlength{\extrarowheight}{-10pt}
\maketitle

\begin{abstract}
We introduce a novel strategy to address the issue of demand estimation in single-item single-period stochastic inventory optimisation problems. Our strategy analytically combines confidence interval analysis and inventory optimisation. We assume that the decision maker is given a set of past demand samples and we employ confidence interval analysis in order to identify a range of candidate order quantities that, with prescribed confidence probability, includes the real optimal order quantity for the underling stochastic demand process with unknown stationary parameter(s). In addition, for each candidate order quantity that is identified, our approach can produce an upper and a lower bound for the associated cost. We apply our novel approach to three demand distribution in the exponential family: binomial, Poisson, and exponential. For two of these distributions we also discuss the extension to the case of unobserved lost sales. Numerical examples are presented in which we show how our approach complements existing frequentist --- e.g. based on maximum likelihood estimators --- or Bayesian strategies.\\
{\bf keywords:} inventory control, newsvendor problem, confidence interval analysis, demand estimation, sampling.
\end{abstract}

\section{Introduction}

We consider the problem of controlling the inventory of a single item with stochastic demand over a single period. This problem is known as the ``newsvendor'' problem \citep{spp98}. Most of the literature on the newsboy problem has focused on the case in which the demand distribution and its parameters are known. However, what happens in practice is that the decision maker must estimate the order quantity from a --- possibly very limited --- set of past demand realisations. This task is often complicated by the fact that unobserved lost sales must be taken into account. 

Existing approaches to the newsvendor under limited historical demand data can be classified in non-parametric and parametric approaches. 
Non-parametric approaches operate without any access to and assumptions on the true demand distributions. Parametric approaches assume that demand realisations come from a given probability distribution --- of from a family of probability distributions --- and make inferences about the parameters of the distribution. When the class of the distribution is known, but its parameters must be estimated from a set of samples, non-parametric approaches may produce conservative results. For this reason several works in the literature investigated parametric approaches to the newsvendor under limited historical demand data; a complete overview on these works will be provided in Section \ref{sec:literature}. Two classical parametric approaches for dealing with the newsvendor under limited historical demand data are the maximum likelihood \citep[see e.g.][]{Scarf1959,fukuda1960,gupta1960} and the Bayesian approach \citep[see e.g.][]{RePEc:eee:ejores:v:98:y:1997:i:3:p:555-562,Hill99}. Both these strategies feature a number of asymptotical properties that guarantee their convergence towards the optimal control strategy. However, a decision maker finds herself rarely in an asymptotic situation, since only few samples are generally available to estimate an order quantity. This means that asymptotic properties often do not hold in practice. Unfortunately, both the maximum likelihood estimation and the Bayesian approaches ignore the uncertainty around the estimated order quantity and its associated expected total cost or profit. \cite{Hayes1969} and, more recently, \cite{citeulike:12594158} discuss how to quantify this uncertainty by using the concept of expected total operating cost (ETOC), which represents the expected one-period cost associated with operating under an estimated inventory policy. By minimising this performance indicator, they identify an optimal biased order quantity that accounts for the uncertainty around the demand parameters estimated from limited historical data. Their approach, however, does not answer a number of fundamental questions. It does not state at what confidence level we can claim this order quantity to be optimal within a given margin of error; nor does it quantifies the probability of incurring an expected total cost substantially different than the estimated one, when such an order quantity is selected. 
\cite{Kevork_10} exploits the sampling distribution of the estimated demand parameters to study the variability of the estimated optimal order quantity and its expected total profit under a normally distributed demand with unknown parameters. The author shows that these two estimators asymptotically converge to normality. Based on this property, asymptotic confidence intervals are derived for the true optimal order quantity and its expected total profit. Unfortunately, these confidence intervals achieve the prescribed confidence level only asymptotically and they represent an approximation when one operates under finite samples. In this work, we make the following contributions to the inventory management literature:
\begin{enumerate}
\item We analytically combine confidence interval analysis and inventory optimisation. By exploiting exact confidence intervals for the parameters of a given distribution, we identify a range of candidate order quantities that, according to a prescribed confidence probability, includes the real optimal order quantity for the underlying stochastic demand process with unknown stationary parameter(s). For each candidate optimal order quantity that is identified, our approach computes an upper and a lower bound for the associated cost. This range covers, once more according to a prescribed confidence probability, the actual cost the decision maker will incur if she selects that particular quantity.
\item To obtain the former result, when the order quantity is fixed, we establish convexity of the newsvendor cost function in the success probability $p$ of a binomial demand (Theorem \ref{th:convexity_gq_binomial}) and in the rate $\lambda$ of a poisson demand (Theorem \ref{th:convexity_gq_poisson}); we also establish that the newsboy cost function is quasi-convex in the expected value $1/\lambda$ of an exponential distribution (Theorem \ref{th:theorem_g}). These results are nontrivial and, to the best of our knowledge, they have not been established before in the literature. 
\item For the binomial and the poisson distribution we demonstrate how to extend the discussion to account for unobserved lost sales.
\end{enumerate}

Our strategy is frequentist in nature and based on the theory of statistical estimation introduced by \cite{Ney37}. In contrast to Bayesian approaches, no prior knowledge is required to perform the computation, which is entirely data driven. In contrast to \citep{Hayes1969,citeulike:12594158} we do not introduce new performance indicators, such as the ETOC, and we build our analysis on existing and well established results from inventory theory, i.e. the expected total cost of a policy; and from statistical analysis, i.e. confidence intervals. Finally, in contrast to \citep{Kevork_10} our results are valid both asymptotically and under finite samples.

If the identified set of candidate optimal order quantities comprises more than a single value, expert assessment or any existing frequentist or Bayesian approach may be employed to select the most promising of these values according to a given performance indicator. By using our approach, the decision maker may then determine --- at a given confidence level and from a limited set of available data --- the exact confidence interval for the expected total cost associated with such a decision, as well as the potential discrepancy between the true optimal decision and the one she selected. For this reason, a further contribution is the following.
\begin{enumerate}
\setcounter{enumi}{3}
\item We effectively complement a number of existing strategies that compute a point estimate of the optimal order quantity and its expected total cost. We demonstrate this fact for the Bayesian approach in \citep{RePEc:eee:ejores:v:98:y:1997:i:3:p:555-562} and for the frequentist approach based on the maximum likelihood estimator of the demand distribution parameter.
\end{enumerate}

\section{Literature survey}\label{sec:literature}

A thorough literature review on the newsvendor problem is presented by \cite{Khouja_2000}. Among other extensions, the author surveyed those dealing with different states of information about demand. As \cite{spp98} point out, in practice demand distribution may be not known. \cite{Khouja_2000} points out that several authors relaxed the assumptions of having a specific distribution with known parameters.  

One of the earliest approaches to dealing with different states of information about demand is the so-called ``maximin approach'', which consists in maximising the worst-case profit. \cite{scarf1958min} applied the maximin approach to the newsvendor problem and assumed that only the mean and the variance of demand are known. Under this assumption, they derived a closed-form expression for the optimal order quantity that maximises the expected total profit against the worst-case demand distribution. \cite{Gallego1993} provided a simplified form of the rule in \citep{scarf1958min}; beside this, they also discussed four extensions: a recourse case, a model including fixed ordering cost, a random yield case and a multi-product case.  This model was extended to account for balking in \cite{Moon1995}. \cite{citeulike:12228279} assumed that the demand distribution belongs to a certain class of probability distribution functions with known mean and standard deviation; the authors' aim is to compute the maximum expected value of distribution information over all possible probability distribution functions with known mean and standard deviation for any order quantity. \cite{citeulike:12232091} pointed out that the maximin objective is, generally speaking, conservative, since it focuses on the worst-case scenario for the demand. The authors therefore suggested adopting a less conservative approach: the ``minimax regret'', in which the firm minimises its maximum cost discrepancy from the optimal decision. Works mentioned so far focused on a newsvendor setting. In contrast, \cite{techreport1970} discuss a minimax multi-period inventory model. 
\cite{citeulike:12602650} discuss a finite horizon inventory model in which the demand distribution is discrete and partially defined by selected moments and/or percentiles. \cite{citeulike:6644349} analyse a distribution free inventory problem for which demand in each period is a random variable over a given support identified by two values: the lower and the upper estimators. A comparable model is found in \cite{citeulike:7362522}. \cite{citeulike:12594160} minimise a coherent risk measure instead of the total cost in an inventory control model and establish an equivalence between risk aversion described as a coherent risk measure and a minimax formulation of the problem.   \cite{citeulike:7682492} discuss an inventory control problem under a demand model described by a given support, mean and standard deviation. They then consider a second-order cone optimisation problem that minimises the expected total cost among all distributions satisfying the demand model. 

All the aforementioned works operate in a {\em distribution free} setting. The decision maker has access to partial information about the demand distribution, e.g. mean, variance, symmetry, unimodality etc, but does not know the class of the demand distribution, e.g. Poisson, normal etc. In practice, it is often the case that the decision maker can only access a set of past observations of the demand out of which a near-optimal inventory target must be estimated. Approaches trying to estimate a near-optimal inventory target from observed realisations of the demands can be classified as {\em non-parametric} or {\em parametric}. 

Non-parametric approaches operate without any access to and assumptions on the true demand distributions. \cite{citeulike:12602630} discuss a non-parametric approach which computes policies based only on observed samples of the demands. The authors determine bounds for the number of samples that are necessary in order to guarantee an arbitrary approximation of the optimal policy defined with respect to the true underlying demand distributions. \cite{citeulike:12602635} discuss an adaptive inventory policy that deals with censored observations, thus effectively relaxing the assumption that past demand data is observable. Other non-parametric approaches based on order statistics were proposed in \citep{citeulike:12602636,citeulike:12602637}; approaches based on bootstrapping techniques were discussed in \citep{citeulike:12602638,citeulike:12602639}. 

Parametric approaches assume that demand realisations come from a given probability distribution and make inferences about the parameters of the distribution. The class of the distribution may be determined, for instance, by selecting the maximum entropy distribution that matches the structure of the demand process at hand (see also the discussion in \cite{citeulike:12232091} p. 190). According to Jaynes' ``principle of maximum entropy'', introduced in \citep{1957PhRv}, subject to known constraints (i.e. testable information), the probability distribution which best represents the current state of knowledge is the one with largest entropy. When the class of the distribution is known but its parameters must be estimated from a set of samples non-parametric approaches such as those discussed so far are not appropriate, since they would produce conservative results.  

According to \cite{RePEc:eee:ejores:v:182:y:2007:i:1:p:256-281} there are two general approaches for dealing with a stochastic decision making environment in which random variables follow known distributions with unknown parameters: the Bayesian and the frequentist. 
In the Bayesian approach  a ``prior'' distribution is selected for the demand distribution parameter(s). This selection may be based on collateral data and/or subjective assessment. Subsequently, a ``posterior'' distribution is derived from the prior distribution by using the demand data that is observed.  This posterior distribution is used to construct the posterior distribution of the demand and to determine the optimal order quantity and objective function value.
In the frequentist approach a parametric demand distribution is empirically selected and point estimates, e.g. maximum likelihood or moment estimators, for the unknown parameters are obtained according to the observed data; these are then used to derive the optimal order quantity and objective function value. 

According to \cite{Kevork_10} another distinction can be made between approaches assuming that demand is fully observed and approaches assuming that demand occurring when the stock level drops to zero is lost and thus not observed. In the latter case, it is necessary to adjust the estimation procedure to account for unobserved demand. \cite{citeulike:12237796} further distinguish works on estimating demand distributions with stockouts in two groups: estimating the initial demand distribution, e.g. \cite{Conrad1976,citeulike:12237796}; and updating the demand distribution parameters \cite{Wecker1978,Bell1981,RePEc:eee:proeco:v:28:y:1992:i:2:p:211-215}. 

Bayesian approaches under fully observed demand were proposed in \cite{Scarf1959,NAV:NAV3800070428,Iglehart1964,Azoury1985,Lovejoy1990,Bradford1990,RePEc:eee:ejores:v:98:y:1997:i:3:p:555-562,Eppen1997,Hill99,citeulike:12228282,DBLP:journals/rda/BensoussanCRS09}. Bayesian approaches under censoring induced by lost sales include \cite{citeulike:12237818,Ding98thecensored,RePEc:eee:ejores:v:182:y:2007:i:1:p:256-281,citeulike:12228284,Lu2008,Mersereau2012}. Bayesian approaches suffer from a number of drawbacks. First, an initial ``belief'' about the unknown parameters must be expressed as a prior distribution of the demand. It is often assumed that this prior distribution is obtained from collateral data and/or subjective assessment. The need of a prior distribution is structural in the Bayesian approach, which relies on the update of the prior distribution to derive the posterior distribution of the demand given the data. When no supporting information is available, ``uninformative'' priors can be used, see e.g. \cite{Hill99}, but these tend to introduce a strong bias expecially under limited available data to perform Bayesian updating. This fact is well known in the life sciences, see e.g. \cite{citeulike:2184775}, but it is often ignored in more theoretical settings. At the end of this work, in Section \ref{sec:dis}, we will provide a practical exemplification of this fact. A second issue that arises with existing Bayesian approaches to the newsvendor problem is that to show that the order quantity derived via the Bayesian approach converges to the optimal order quantity one has to consider an infinite set of samples, see e.g. \cite{DBLP:journals/rda/BensoussanCRS09}. However, in practice it is often the case that available information is very limited. Unfortunately, Bayesian approaches can be shown to be often biased under small sample sets, especially due to the fact that the choice of the prior may strongly influence the order quantity obtained.

Two early frequentist approaches are \cite{Nahmias_94} and \cite{Agrawal_96}. \cite{Nahmias_94} analyzed the censored demand case under a normally distributed demand. \cite{Agrawal_96} considered the estimation of a negative binomial demand under censoring induced by lost sales. However, these two studies consider the stock level as given and thus do not address the associated optimization problem of finding the optimum stock level. More recently, \cite{citeulike:12594157} introduced the ``operational statistics'' framework, in which an optimal order quantity, rather than the parameters of the distribution, is directly estimated from the data. The authors consider the case in which it is known that the demand distribution function belongs to a parameterised family of distribution functions. In contrast to the Bayesian approach, they do not assume any prior knowledge on the parameter values. They demonstrate that by combining demand parameter estimation and inventory-target optimisation a higher expected total profit can be achieved with respect to traditional approaches that separate estimation and optimisation. 
\cite{citeulike:12594156} integrate distribution fitting and robust optimisation by identifying a set of demand distributions that fit historical data according a given criterion; they then characterise an optimal policy that minimises the maximum expected total cost against such set of demand distributions.

There is an important limitation that is common to all approaches surveyed so far. Consider a possibly very limited set of past demand observations; strategies based on frequentist analysis, e.g. maximum likelihood estimators and distribution fitting; or on Bayesian analysis, e.g. \cite{Hill99}, would analyze these data and provide a single most-promising order quantity and an estimated cost associated with it. However, given the available data, they do not answer a number of fundamental questions: we do now know at what confidence level we can claim the quantity selected to be optimal within a given margin of error; and we also do not know the probability of incurring a cost substantially higher than the estimated one, when such an order quantity is selected.

To the best of our knowledge \cite{Kevork_10} was the first to exploit the sampling distribution of the estimated demand parameters in order to study the variability of the estimates for the optimal order quantity and associated expected total profit. The author adopts a frequentist approach in which demand is fully observed in each period. By incorporating maximum likelihood estimators for mean and variance of demand into expressions that determine the optimal order quantity and associated expected total profit, the author develops estimators for these latter two variables. These estimators are shown to be consistent and to asymptotically converge to normality. Based on these properties, the author derives confidence intervals for the true optimal order quantity and associated expected total profit. Unfortunately, these estimators are biased in finite samples and the associated confidence intervals achieve the prescribed confidence level only asymptotically.

As pointed out in \citep{citeulike:12594158}, the inventory manager finds rarely herself in an asymptotic situation, since an inventory target must be typically estimated from a small sample size. To quantify the uncertainty about distribution parameter estimates and thus about the estimated order quantity, \cite{citeulike:12594158} adopt the ETOC, originally introduced in \citep{Hayes1969}, which we recall represents the expected one-period cost associated with operating under an estimated inventory policy. Originally, \cite{Hayes1969} discussed applications of ETOC to exponentially and normally distributed demands. More specifically, they identified the optimal biased order quantity that minimises the ETOC in presence of limited historical demand data. This was one of the first works blending statistical estimation with inventory optimisation. \cite{citeulike:12594158} extended this analysis to a parameterised family of distributions --- the Johnson translation system --- that has the ability to match any finite first four moments of a random variable and to capture a broad range of distributional shapes. The authors quantify the inaccuracy in the order quantity estimation by using the expected value of perfect information about the sampling distribution of the demand parameters for the estimated order quantity. By using this interpretation of the ETOC, they seek an order quantity that minimises the ETOC within a class of inventory target-estimators implied by the Johnson translation system. Despite its ability to quantify the inaccuracy in the inventory-target estimation as a function of the length of the historical data via the ETOC, the approach in \citep{citeulike:12594158} does not identify a confidence interval that, with prescribed confidence probability, includes the real optimal order quantity for the underlying stochastic demand process with unknown parameter(s); neither it is able to produce a confidence interval for the expected total cost associated with ordering decisions in this interval.

\section{The newsvendor problem}\label{sec:new}

In this section, we shall summarize the key features of the newsvendor problem.
For more details, the reader may refer to \cite{spp98}.
Consider a one-period random demand $d$ with mean $\mu$ and variance $\sigma^2$.
Let $h$ be the unit overage cost, paid for each item left in stock after demand is realized,
and let $p$ be the unit underage cost, paid for each unit of unmet demand; we assume 
$p,h>0$. 
Let $g(x)=hx^+ + p x^-$, where $x^+=\max(x,0)$ and $x^-=-\min(x,0)$. The expected
total cost can be written as $G(Q)=E[g(Q-d)]$, where $E[\cdot]$ denotes the expected value. 

It is a well-known fact that, for a continuous demand $d$, 
\begin{equation}\label{eq:critical_fractile_1}
F(Q)=\Pr\{d\leq Q\}=\frac{p}{p+h}=\beta.
\end{equation}
If $F$ is continuous there is at least one $Q$ satisfying Eq. \ref{eq:critical_fractile_1}, that is 
\[Q^*=\inf\{Q\geq0:F(Q)=\beta\}.\]
For $F$ strictly increasing, 
there is a unique 
optimal solution given by 
\begin{equation}\label{eq:critical_fractile}
Q^*=F^{-1}(\beta).
\end{equation}

In practice, the probability distribution of the random demand $d$ often has finite support 
over the set $\mathbb{N}_0=\{0,1,2,\dots\}$. In this case it is useful to work with the forward difference
$\Delta G(Q)= G(Q+1)-G(Q)$, $Q\in \mathbb{N}_0$. 
It is easy to see that 
\[\Delta G(Q)=h-(h+p) \Pr\{d>j\}\]
is non-decreasing in $Q$, and that $\lim_{Q\rightarrow \infty}\Delta G(Q)=h>0$, so an optimal solution is given by $Q^*=\min\{Q\in\mathbb{N}_0:\Delta G(Q)\geq 0\}$ or equivalently
\begin{equation}\label{eq:critical_fractile_disc}
Q^*=\min\{Q\in\mathbb{N}_0:F(Q)\geq \beta\}.
\end{equation}

\subsection{A frequentist and a Bayesian approach}\label{sec:heuristics}

Let us now consider the situation in which the decision maker knows the
type of the random demand distribution (e.g. binomial), but in which he
does not know the actual values of some or all the (stationary) parameters of such a
distribution. Nonetheless, the decision maker is given a set of $M$ past 
realizations of the demand. From these realizations
he has to infer the optimal order quantity and, possibly, he has to estimate the cost 
associated with the optimal $Q^*$ he has selected. 

In what follows, we detail the functioning of a frequentist approach, i.e. the maximum likelihood approach, and of a Bayesian approach from the literature \cite{RePEc:eee:ejores:v:98:y:1997:i:3:p:555-562}. In the rest of this work we will make ample use of these two approaches, especially to discuss how our approach can be used to complement the results obtained by frequentist or Bayesian approaches.  For the sake of brevity, in this work we will focus only on these two strategies from the literature. However, this choice is made without loss of generality. Our approach may in fact also complement any of the other frequentist or Bayesian approaches previously surveyed.

\subsubsection{Maximum likelihood approach}

A commonly adopted heuristic strategy for order quantity selection
under sampled demand information consists in computing, from the
available sample set,  a point
estimate for the unknown demand distribution parameter(s). 
This may be done by using the maximum likelihood estimator \citep{citeulike:7982953}, thus obtaining the so-called maximum likelihood policy \citep[see e.g.][]{Scarf1959,fukuda1960,gupta1960}, or
the method of moments \citep{RePEc:eee:ecochp:4-36}. For instance, assume that the available
sample set comprises $M$ observed past demand data,
$d_1,\ldots,d_M$, and that the demand is assumed to follow
a binomial distribution. The binomial distribution comprises
two parameters: the number of trials $N$ and the success probability $q$.
In the context of the newsvendor problem, $N$ may be interpreted as 
the non-variant number of customers that enter the shop in a given 
day $i\in\{1,\ldots,M\}$, and $q$ may be interpreted as the probability that a customer
makes a purchase. Then, by assuming that $N$ is known, the maximum likelihood 
estimator for parameter $q$ is $\widehat{q}=\sum_{i=1}^M d_i/(M N)$.
After computing $\widehat{q}$, the decision maker employs 
the random variable $\text{bin}(N;\widehat{q})$ in place of the actual
unknown demand distribution in Eq. \ref{eq:critical_fractile_disc} 
to compute the estimated optimal order quantity $\widehat{Q}^*$ and 
expected total cost $G(\widehat{Q}^*)$.

\subsubsection{Hill's Bayesian approach}

A Bayesian approach to the Newsvendor problem under partial information is
presented by \cite{RePEc:eee:ejores:v:98:y:1997:i:3:p:555-562}. 
Hill's approach considers a ``prior'' distribution, based on collateral data and/or
subjective assessment, for the unknown parameter(s) of the demand distribution.
As new data are observed, the prior distribution is updated and a ``posterior'' 
distribution is generated. Hill then uses the posterior distribution of the unknown
parameter(s) for estimating the posterior distribution of the demand. Finally, the
posterior distribution of the demand is used to estimate the order quantity 
that optimizes the Newsvendor cost function.
More formally, we recall that Bayes' theorem tells us that 
\[\Pr\{a_j|b\}=\frac{\Pr\{b|a_j\}\cdot \Pr\{a_j\}}{\sum_{i=1}^k\Pr\{b|a_i\}\cdot\Pr\{a_i\}},\]
where $\{a_1,\ldots,a_k\}$ is a partition of the sample space and
$b$ is an observed event. 
Bayes actually discusses also the natural extension of the above theorem
when $a$ is continuous and $b$ is discrete or continuous,
\[f(a|b)=\frac{f(b|a)\cdot g(a)}{\int f(b|u)\cdot g(u)\mathrm{d}u}.\]
The denominator is, of course, independent of $a$, therefore
$f(a|b)\propto f(b|a)\cdot g(a)$.
In the context of the Newsvendor problem, $a$ represents the unknown
parameter of the demand distribution, $b$ represents the actual set
of observed demand samples. According to Hill, the prior distribution of $a$, $g(a)$,
describes an estimate of the likely value that $a$ might take, without considering
the observed samples. This estimation is based on subjective 
assessment and/or collateral data. $f(b|a)$, also known as the likelihood, represents
the probability of observing a set of samples $b$ given $a$. The posterior
distribution of $a$, $f(a|b)$, is an updated estimate of the values $a$ is likely to take
based on the prior distribution and the observed information.
Hill uses an uninformative, also known as objective, prior to express ``an initial state of 
complete ignorance of the likely values that the parameter $a$ might take.'' By employing
the conjugate prior for the particular distribution under analysis, he constructs the
posterior distribution for the Newsvendor demand as follows,
\[f(d|b)=\int f(d|a)f(a|b)\mathrm{d}b,\]
where the integral is computed over all permitted values of $a$.
The Bayesian approach proposed by Hill then consists in using this posterior distribution
in place of the unknown true distribution for the demand in Eq. \ref{eq:critical_fractile}.
This immediately produces a candidate order quantity $\widehat{Q}^*$.

\section{Binomial demand}\label{sec:bin}

Consider a discrete random variable that follows a Bernoulli distribution.
Accordingly, such a variable may produce only two possible outcomes, i.e. ``yes''
and ``no'', with probability $q$ and $1-q$, respectively. 
This class of random variables is particularly useful in representing
so called ``Bernoulli trials'', which are experiments that can have one of two possible
outcomes. These  events can be phrased into ``yes or no'' questions, such as ``did the 
customer purchase the newspaper?''

Consider the following situation: a Newsvendor has a pool of $N$ customers that come
every day to the stand. Each customer may buy a newspaper with probability 
$q$. It is a well-known fact that any experiment comprising a sequence of $N$ 
Bernoulli trials, each having the same
``yes''  (respectively, ``no'') probability $q$ (respectively, $1-q$), can be represented
by a random variable $\text{bin}(N; q)$ that follows a binomial distribution \citep{j61} with
probability mass function
\[\Pr\{\text{bin}(N; q)=k\}= \binom{N}{k}q^k(1-q)^{N-k},\]
where $k=0,\ldots,N$.

According to our previous discussion it is fairly 
easy to find the optimal order quantity $Q^*$ for a given random demand $\text{bin}(N; q)$. 
We shall now give a running example. 
Consider a random demand $\text{bin}(50,0.5)$. Let
$h=1$ and $p=3$, therefore $\beta=0.75$. From Eq. \ref{eq:critical_fractile_disc} 
we compute $Q^*=27$ and
the respective expected total cost $G(Q^*)=4.4946$.

Let us now consider the situation in which the parameter $q$ is not known. 
The decision maker is given a set of $M$ past realizations of the demand. From these realizations
he has to infer the optimal order quantity and, possibly, he has to estimate the cost 
associated with the optimal $Q^*$ he has selected.

Since we operate under partial information it may not be possible to uniquely
determine ``the'' optimal order quantity and the exact cost associated with it. Therefore,
we argue that a possible approach consists in determining a range of ``candidate'' optimal
order quantities and upper and lower bounds for the cost associated with these quantities.
This range will contain the real optimum according to a prescribed confidence probability.

\subsection{Confidence intervals for the binomial distribution}

Confidence interval analysis \citep{Ney37,ney41} is a well established technique in statistics 
for computing, from a given set of experimental results, a range of values 
that, with a certain confidence level (or confidence probability), 
will cover the actual value of a parameter that is being estimated.
Several techniques \citep[etc.]{cloppears34,Garwood36,trivedi01} for building
confidence intervals for a given sample set have been proposed.

Approximate techniques for building confidence intervals \citep[see e.g.][]{agresticloull98}
become relevant because, especially with small sample sizes, an exact
confidence interval may be unnecessarily conservative. In this work, we focus on the exact interval and we leave the 
analysis of the benefits brought by the use of approximate intervals as future research.

A method to compute ``exact confidence intervals'' for the 
binomial distribution has been introduced by \cite{cloppears34}.
This method uses the binomial cumulative distribution function in order to build the interval
from the data observed. The Clopper-Pearson interval can be written as
$\left(q_{lb},q_{ub}\right)$, 
where
\[
\begin{array}{l}
q_{lb}=\min\{q| \Pr\{\text{bin}(N;q)\geq X\}\geq (1-\alpha)/2\},\\
q_{ub}=\max\{q| \Pr\{\text{bin}(N;q)\leq X\}\geq (1-\alpha)/2\},
\end{array}
\]
$X$ is the number of successes (or ``yes'' events) observed in the sample and $\alpha$ is the confidence probability. Note that we assume $q_{lb}=0$ when $X=0$ and $q_{ub}=N$ when $X=N$. As discussed
by \cite{agresticloull98}, this interval can be also expressed using quantiles from the beta distribution.
More specifically, the lower endpoint is the $(1-\alpha)/2$-quantile of a beta distribution with shape parameters 
$X$ and $N-X+1$, and the upper endpoint is the $(1+\alpha)/2$-quantile of a beta distribution with shape parameters $X+1$ and $N-X$.
Furthermore, the beta distribution is, in turn, related to the F-distribution so a third formulation of the Clopper-Pearson interval, also discussed by \cite{agresticloull98}, uses quantiles from the F distribution. 

It is intuitively clear that the ``quality'' of a given confidence interval is directly 
related to its size. The smaller the interval, the better the estimate.
In general, confidence intervals that have symmetric tails (i.e. 
with associated probability $(1-\alpha)/2$) are not the smallest possible ones. 
A large literature exists on the topic of determining the smallest possible intervals
for a given parameter/distribution combination \citep[see e.g.][]{citeulike:7874841}. 
The discussion that follows is independent of the particular interval adopted. 
For the sake of simplicity, we will adopt intervals having symmetric tails. 

\subsection{Solution method employing statistical estimation based on classical theory of probability}\label{sec:solution_method}

We shall now employ the Clopper-Pearson interval for computing an upper and a lower bound
for the optimal order quantity $Q^*$ in a Newsvendor problem under partial information. 
The confidence interval for the unknown parameter $q$ of the binomial demand $\text{bin}(N;q)$ is simply $\left(q_{lb},q_{ub}\right)$
where
\[
\begin{array}{l}
q_{lb}=\min\{q| \Pr\{\text{bin}(MN;q)\geq X\}\geq (1-\alpha)/2\},\\
q_{ub}=\max\{q| \Pr\{\text{bin}(MN;q)\leq X\}\geq (1-\alpha)/2\},
\end{array}
\]
and $X=\sum_{i=1}^M d_i$. Let $Q^*_{lb}$ be the optimal order quantity for the Newsvendor problem under a $\text{bin}(N,q_{lb})$ demand and $Q^*_{ub}$ be the optimal order quantity for the Newsvendor problem under a $\text{bin}(N,q_{ub})$ demand. Since $\Delta G(Q)$ is non-decreasing in $Q$,
according to the available information with confidence probability $\alpha$ the optimal order quantity $Q^*$ is a member of the set $\{Q^*_{lb},\ldots,Q^*_{ub}\}$. 

We shall now compute upper ($c_{ub}$) and lower ($c_{lb}$) bounds for
the cost associated with a solution that sets the order quantity to a
value in the set $\{Q^*_{lb},\ldots,Q^*_{ub}\}$.  Let us write the
cost associated with an order quantity $Q$,
\[G(Q)=h\sum_{i=0}^{Q} \Pr\{\text{bin}(N;q)=i\} (Q-i) + p\sum_{i= Q}^N \Pr\{\text{bin}(N;q)=i\} (i-Q).\]
Then, consider the function 
\begin{equation}\label{eq:G(q)}
G_Q(q)=h\sum_{i=0}^{Q} \Pr\{\text{bin}(N;q)=i\} (Q-i) + p\sum_{i= Q}^N \Pr\{\text{bin}(N;q)=i\} (i-Q),
\end{equation}
in which the order quantity $Q$ is fixed and in which we vary the
``success'' probability $q\in(0,1)$.  It can be proved that $G_Q(q)$ is
convex in the continuous parameter $q$; this is trivially true when
$Q=N$.  The proof for $0\leq Q<N$ is given in Appendix \ref{sec:appendix_I} (Theorem \ref{th:convexity_gq_binomial}).

Although it is possible to prove that $G_Q(q)$ is convex in $q$, there is no closed form expression for finding the $q^*$
that minimizes this function. Nevertheless, due to its convexity in $q$, it is clearly 
possible to use convex optimization approaches to find 
the $q^*$ that minimizes or maximizes this function over a given interval.

Let us consider the confidence interval $\left(q_{lb},q_{ub}\right)$ 
for the parameter $q$ of the binomial demand. 
For a given order quantity $Q$, 
consider the value \[q^*_{Q,\min}=\argmin_{q\in(q_{lb},q_{ub})}G_Q(q)~~~~~~~~~~~~~~
\left(q^*_{Q,\max}=\argmax_{q\in(q_{lb},q_{ub})}G_Q(q)\right)\] that minimizes (maximizes) 
$G_Q(q)$ for $q\in(q_{lb},q_{ub})$. With confidence probability $\alpha$,
$G_Q(q_{Q,\min}^*)$ and $G_Q(q_{Q,\max}^*)$ represent a lower and an upper bound,
respectively, for the cost associated with $Q$. 

By recalling that the optimal order quantity $Q^*$ is, with confidence probability
$\alpha$, a member of the set $\{Q^*_{lb},\ldots,Q^*_{ub}\}$, it is easy to compute
upper ($c^*_{ub}$) and lower ($c^*_{lb}$) bounds for the cost that a manager 
will face, with confidence probability $\alpha$, whatever order quantity he 
chooses in the candidate set $\{Q^*_{lb},\ldots,Q^*_{ub}\}$. 
The lower bound is
\[c^*_{lb}=\min_{Q\in\{Q^*_{lb},\ldots,Q^*_{ub}\}}{G_Q(q^*_{Q,\min})}\]
and the upper bound is
\[c^*_{ub}=\max_{Q\in\{Q^*_{lb},\ldots,Q^*_{ub}\}}{G_Q(q^*_{Q,\max})}.\]

It should be emphasized that, when the confidence interval $\left(q_{lb},q_{ub}\right)$ 
covers the real parameter $q$ of the binomial demand we are estimating, then the set 
$\{Q^*_{lb},\ldots,Q^*_{ub}\}$ comprises the optimal order quantity $Q^*$
and the interval $(c^*_{lb},c^*_{ub})$ comprises the real cost associated with
every possible order quantity in $\{Q^*_{lb},\ldots,Q^*_{ub}\}$. Given the way confidence
interval $\left(q_{lb},q_{ub}\right)$ is constructed, it is guaranteed that this 
happens with probability $\alpha$. 

Of course, by increasing the number $M$ of past observations, we can decrease
the size of confidence interval $\left(q_{lb},q_{ub}\right)$. As a direct consequence,
the cardinality of the set $\{Q^*_{lb},\ldots,Q^*_{ub}\}$ decreases. In the ideal case,
this set comprises a single candidate order quantity $Q^*$ that with confidence
probability $\alpha$ represents an optimal solution to the problem and has a 
cost comprised in the interval $(G_{Q^*}(q^*_{Q^*,\min}),G_{Q^*}(q^*_{Q^*,\max}))$. 

Finally, consider the case in which unobserved lost sales occurred and the 
$M$ observed past demand data, $d_1,\ldots,d_M$, only reflect the number of customers that
purchased an item when the inventory was positive. The analysis discussed
above can still be applied provided that the confidence interval 
for the unknown parameter $q$ of the $\text{bin}(N;q)$ demand is computed as
\[
\begin{array}{l}
q_{lb}=\min\{q| \Pr\{\text{bin}(\sum_{j=1}^M \widehat{N}_j;q)\geq X\}\geq (1-\alpha)/2\},\\
q_{ub}=\max\{q| \Pr\{\text{bin}(\sum_{j=1}^M \widehat{N}_j;q)\leq X\}\geq (1-\alpha)/2\},
\end{array}
\]
where $\widehat{N}_j$ is the total number of customers that entered the shop in 
day $j$ --- for which a demand sample $d_j$ is available --- while the inventory
was positive.

\subsection{Algorithm}

The procedure to compute, under the prescribed confidence probability $\alpha$, 
a candidate set $\mathcal{Q}$ of order quantities and 
upper $(c^*_{ub})$ and lower $(c^*_{lb})$ bounds for the cost
a manager faces when he selects one of these quantities is presented in Algorithm
\ref{newsvendor_binomial_alg}. 
\begin{figure}[h!]
\vskip 1pc
  \begin{algorithm}[H]
    
    \SetKwInOut{Input}{input}
    \SetKwInOut{Output}{output}
    \Input{$M$ past demand realizations $d_i$, $i=1,\ldots,M$;\\
    the number of customers per day: $N$;\\
    the holding cost: $h$;\\
    the penalty cost: $p$;\\
    the confidence probability: $\alpha$.}
    \Output{the set $\mathcal{Q}$ of candidate order quantities;\\
    the interval $(c^*_{lb},c^*_{ub})$ for the estimated cost.}
    \BlankLine
    \Begin{
    	$a\leftarrow \sum_{i=1,\ldots,M}d_i$\;
	$b\leftarrow M N-\sum_{i=1,\ldots,M}d_i$\;
         $q_{ub}\leftarrow$\texttt{InverseCDF[BetaDistribution($a+1,b$),$(1+\alpha)/2$]}\;
         $q_{lb}\leftarrow$\texttt{InverseCDF[BetaDistribution($a,b+1$),$(1-\alpha)/2$]}\;
         $\beta\leftarrow p/(p+h)$\;
         $Q^*_{ub}\leftarrow$\texttt{InverseCDF[BinomialDistribution($N,q_{ub}$),$\beta$]}\;
         $Q^*_{lb}\leftarrow$\texttt{InverseCDF[BinomialDistribution($N,q_{lb}$),$\beta$]}\;
         $\mathcal{Q}\leftarrow \{Q^*_{lb},\ldots,Q^*_{ub}\}$\;
    	$c^*_{ub}\leftarrow -\infty$\;
	$c^*_{lb}\leftarrow \infty$\;
	\For{each $Q\in\mathcal{Q}$}{
		$q^*_{Q,\max}\leftarrow \argmax_{q\in(q_{lb},q_{ub})}G_Q(q)$\;
		$q^*_{Q,\min}\leftarrow \argmin_{q\in(q_{lb},q_{ub})}G_Q(q)$\;
		$c^*_{ub}\leftarrow\max(c^*_{ub},G_Q(q^*_{Q,\max}))$\;
		$c^*_{lb}\leftarrow\min(c^*_{lb},G_Q(q^*_{Q,\min}))$\;
	}  
    }
    \caption{Newsvendor under incomplete information: binomial demand. \label{newsvendor_binomial_alg}}
  \end{algorithm}
\end{figure}
The code initially computes Clopper-Pearson interval ($q_{lb},q_{ub}$) by exploiting the
relationship between the binomial distribution and the Beta distribution \citep{ehp00} --- \texttt{InverseCDF} denotes the inverse cumulative distribution function. Then
it computes the critical fractile $\beta$ and the upper and lower bound for the set $\mathcal{Q}$
of candidate order quantities. Finally, it iterates through the elements of this set to 
compute the upper ($c^*_{ub}$) and lower bound ($c^*_{lb}$) for the estimated cost
associated with these candidate order quantities.

In general, the set $\mathcal{Q}=\{Q^*_{lb},\ldots,Q^*_{ub}\}$
may comprise a significant number of elements, especially 
if a very limited number of samples is available. 
A decision maker may then employ one of the strategies discussed 
in Section \ref{sec:heuristics} in order to determine the most promising quantity in this set. 

\def\registered{{\ooalign{\hfil\raise .00ex\hbox{\scriptsize R}\hfil\crcr\mathhexbox20D}}}

\subsection{Example}\label{sec:example_1}

We consider a simple example involving the Newsvendor problem under binomial demand. 
Assume that, in our problem, $h=1$, $p=3$, and the demand follows a $\text{bin}(50,q)$ distribution, 
in which parameter $q$ is unknown. We are
given $10$ samples for the demand, which we may use to determine the optimal order
quantity $Q^*$. The samples are
$\{28, 28, 24, 27, 25, 26, 28, 28, 23, 27\}$.
The real value for parameter $q$, which is used to generate the 10 samples
is 0.5. Accordingly, the optimal order 
quantity $Q^*$ is equal to 27 and provides a cost equal to 4.4946. 

We consider $\alpha=0.9$. 
By using Algorithm \ref{newsvendor_binomial_alg} we compute the 
set of candidate order quantities $\mathcal{Q}=\{27,28,\ldots,31\}$ and
the confidence interval for the estimated cost $(c^*_{lb},c^*_{ub})=(4.4268,7.2205).$
Among the candidate quantities in $\mathcal{Q}$, both the strategies presented in Section \ref{sec:heuristics} 
identify $\widehat{Q}^*=29$ as the candidate optimal quantity. By using the approach discussed
in Section \ref{sec:solution_method}, we compute the $\alpha$ confidence interval for the estimated cost, 
which is \[\left(G_{\widehat{Q}^*}(q^*_{\widehat{Q}^*,\min}),
G_{\widehat{Q}^*}(q^*_{\widehat{Q}^*,\max})\right)=(4.4487,4.9528).\]

Clearly, the information on the minimum and maximum cost associated with
each order quantity in $\mathcal{Q}$ lets the decision maker perform a more
educated choice. For instance, if a manager is not a risk-taker, he may decide 
select the order quantity  $\widehat{Q}^*$, for which 
the $\alpha$ confidence interval for the estimated cost has the lowest possible upper bound
$G_{\widehat{Q}^*}(q^*_{\widehat{Q}^*,\max})$. In the above example,
this is still 29, but in general it may be a different order quantity.

Less conservative, but approximate, confidence intervals may be
obtained by replacing the Clopper-Pearson \citep{cloppears34} interval with the Agresti-Coull
\citep{agresticloull98} interval for the binomial parameter. The maximum likelihood and the Bayesian
approach do not employ confidence intervals for selecting the candidate order
quantity. Therefore they are not affected by this choice and the selected order quantity 
remains $\widehat{Q}^*=29$. The $\alpha$ confidence interval for the estimated cost 
associated with $\widehat{Q}^*=29$ is
$\left(G_{\widehat{Q}^*}(q^*_{\widehat{Q}^*,\min}),G_{\widehat{Q}^*}
(q^*_{\widehat{Q}^*,\max})\right)=(4.4487, 4.9155)$. This interval is 7.3\% 
smaller than that produced by using the Clopper-Pearson interval.

\section{Poisson demand}\label{sec:poi}

In many practical contexts, a random demand distributed according to a Poisson
law may become relevant. A random demand $\text{Poisson}(\lambda)$ is said to be distributed according to a 
Poisson law with rate parameter $\lambda>0$, if its probability mass function is
\[\Pr\{d=k\}=e^{-\lambda}\frac{\lambda^k}{k!},~~~~~k=0,1,2,\ldots,\infty.\]
The Poisson distribution is the limiting distribution of the binomial distribution when
$N$ is large and $q$ is small. In this case, the parameters of the two distributions are linked by
the relationship $\lambda=qN$. We recall that the expected value of $d$ is $\lambda$
and that the standard deviation of $d$ is $\sqrt{\lambda}$.

By using Eq. \ref{eq:critical_fractile_disc},  we easily obtain the optimal order quantity $Q$ for a 
given demand $d$. We shall give an example. Consider a demand $d$ that follows
a $\text{Poisson}(50)$ distribution. Let $h=1$ and $p=3$, therefore $\beta=0.75$.
The optimal order quantity is $Q^*=55$. Furthermore, by noting that 
\[G(Q)=h(Q-\lambda)+(h+p)\sum_{i=Q}^{\infty}(1-\Pr\{\text{Poisson}(\lambda)\leq i\}),\]
the optimal cost is $G(Q^*)=9.1222$.

We shall now consider, also in this case, the situation in which the parameter $\lambda$ is not known. 
Instead, the decision maker is given a set of $M$ past realizations of $d$. As in the previous case,
from these realizations he has to infer the range of ``candidate'' optimal order quantities and 
upper and lower bounds for the cost associated with these quantities. This range will contain 
the real optimum according to a prescribed confidence probability.

\subsection{Confidence intervals for the Poisson distribution}

As in the previous case, we discuss the exact confidence interval that can be used
to estimate the rate parameter $\lambda$ of the Poisson distribution. This confidence
interval was proposed by \cite{Garwood36} and takes the following form. Consider
a set of $M$ samples $d_i$ drawn from a random demand $d$ that is
distributed according to a Poisson law with unknown parameter $\lambda$. 
We rewrite $\bar{d}=\sum_{i=0}^M d_i$. According to \cite{Garwood36}, the 
confidence interval for $\lambda$ is $(\lambda_{lb},\lambda_{ub})$, where
\[
\begin{array}{l}
\lambda_{lb}=\min\{\lambda| \Pr\{\text{Poisson}(M\lambda)\geq \bar{d}\}\geq (1-\alpha)/2\},\\
\lambda_{ub}=\max\{\lambda| \Pr\{\text{Poisson}(M\lambda)\leq \bar{d}\}\geq (1-\alpha)/2\}.
\end{array}
\]
This interval can be expressed in terms of the chi-square distribution, as shown by \cite{Garwood36}. 
Let $\chi^2_{n}$ denote the chi-square distribution 
with $n$ degrees of freedom, and $G^{-1}(\chi^2_{n},\cdot)$ denote the 
inverse cumulative distribution function of $\chi^2_{n}$. We can write
\[\lambda_{lb}=\frac{G^{-1}(\chi^2_{2\bar{d}},(1-\alpha)/2)}{2M},\]
\[\lambda_{ub}=\frac{G^{-1}(\chi^2_{2\bar{d}+2},(1+\alpha)/2)}{2M}.\]
Furthermore, it is possible to express this interval using quantiles from the gamma distribution \citep{Swift09}. More specifically, the lower endpoint is the $(1-\alpha)/2$-quantile of a gamma distribution with shape parameter $\bar{d}$ and scale parameter $1/M$, and the upper endpoint is the $(1+\alpha)/2$-quantile of a gamma distribution with shape parameter $\bar{d}+1$ and scale parameter $1/M$. Swift lists a number of existing approaches for building 
approximate intervals that are less conservative than Garwood's one and he also
suggests strategies to shorten Garwood's interval by choosing suitable asymmetric tails \citep{Swift09}.

\subsection{Solution method employing statistical estimation based on classical theory of probability}\label{sec:solution_method_poisson}

The method for computing an upper and a lower bound
for the optimal order quantity $Q^*$ in a Newsvendor problem under Poisson demand and partial information
on parameter $\lambda$ can be carried out in a similar fashion to the binomial case given in Section \ref{sec:solution_method}. 
Consider Garwood's confidence interval $\left(\lambda_{lb},\lambda_{ub}\right)$ for the unknown parameter $\lambda$ of the Poisson demand. Let $Q^*_{lb}$ be the optimal order quantity for the Newsvendor problem under a $\text{Poisson}(\lambda_{lb})$ demand and $Q^*_{ub}$ be the optimal order quantity for the Newsvendor problem under a $\text{Poisson}(\lambda_{ub})$ demand. With confidence probability $\alpha$ the optimal order quantity $Q^*$ is a member of the set $\{Q^*_{lb},\ldots,Q^*_{ub}\}$. 

Consider the cost associated with an order quantity $Q$,
\[
G(Q)=h\sum_{i=0}^{Q} \Pr\{\text{Poisson}(\lambda)=i\} (Q-i) + p\sum_{i= Q}^{\infty} \Pr\{\text{Poisson}(\lambda)=i\} (i-Q).
\]
Also in this case we can prove that $G_Q(\lambda)$ 
\begin{equation}\label{eq:G(lambda)}
G_Q(\lambda)=h\sum_{i=0}^{Q} \Pr\{\text{Poisson}(\lambda)=i\} (Q-i) + p\sum_{i= Q}^{\infty} \Pr\{\text{Poisson}(\lambda)=i\} (i-Q),
\end{equation}
is convex in $\lambda$. The proof is given in Appendix \ref{sec:appendix_II} (Theorem \ref{th:convexity_gq_poisson}).
Therefore upper ($c_{ub}$) and lower ($c_{lb}$) bounds for the cost associated with a solution that sets the order quantity to a value in the set $\{Q^*_{lb},\ldots,Q^*_{ub}\}$ can be easily obtained by using convex optimization approaches to find the $\lambda^*$ that minimizes or maximizes this function over a given interval.  

Also in this case, consider the case in which unobserved lost sales occurred and the 
$M$ observed past demand data, $d_1,\ldots,d_M$, only reflect the number of customers that
purchased an item when the inventory was positive. The analysis discussed
above can still be applied provided that the confidence interval 
for the unknown parameter $\lambda$ of the $\text{Poisson}(\lambda)$ demand is computed as
\[
\begin{array}{l}
\lambda_{lb}=\min\{\lambda| \Pr\{\text{Poisson}(\widehat{M}\lambda)\geq \bar{d}\}\geq (1-\alpha)/2\},\\
\lambda_{ub}=\max\{\lambda| \Pr\{\text{Poisson}(\widehat{M}\lambda)\leq \bar{d}\}\geq (1-\alpha)/2\}.
\end{array}
\]
where $\widehat{M}=\sum_{j=1}^M T_j$, and $T_j\in(0,1)$ denotes the fraction of time in day $j$ 
--- for which a demand sample $d_j$ is available --- during which the inventory was positive.

\subsection{Algorithm}

The computational procedure for Poisson demand is presented in Algorithm
\ref{newsvendor_poisson_alg}. 
\begin{figure}[h!]
\vskip 1pc
  \begin{algorithm}[H]
    
    \SetKwInOut{Input}{input}
    \SetKwInOut{Output}{output}
    \Input{$M$ past demand realizations $d_i$, $i=1,\ldots,M$;\\
    the holding cost: $h$;\\
    the penalty cost: $p$;\\
    the confidence probability: $\alpha$.}
    \Output{the set $\mathcal{Q}$ of candidate order quantities;\\
    the interval $(c^*_{lb},c^*_{ub})$ for the estimated cost.}
    \BlankLine
    \Begin{
    	$a\leftarrow \sum_{i=1,\ldots,M}d_i$\;
	$b\leftarrow M$\;
         $\lambda_{ub}\leftarrow$\texttt{InverseCDF[GammaDistribution($a+1,1/b$),$(1+\alpha)/2$]}\;
         $\lambda_{lb}\leftarrow$\texttt{InverseCDF[GammaDistribution($a,1/b$),$(1-\alpha)/2$]}\;
         $\beta\leftarrow p/(p+h)$\;
         $Q^*_{ub}\leftarrow$\texttt{InverseCDF[PoissonDistribution($\lambda_{ub}$),$\beta$]}\;
         $Q^*_{lb}\leftarrow$\texttt{InverseCDF[PoissonDistribution($\lambda_{lb}$),$\beta$]}\;
         $\mathcal{Q}\leftarrow \{Q^*_{lb},\ldots,Q^*_{ub}\}$\;
    	$c^*_{ub}\leftarrow -\infty$\;
	$c^*_{lb}\leftarrow \infty$\;
	\For{each $Q\in\mathcal{Q}$}{
		$\lambda^*_{Q,\max}\leftarrow \argmax_{\lambda\in(\lambda_{lb},\lambda_{ub})}G_Q(\lambda)$\;
		$\lambda^*_{Q,\min}\leftarrow \argmin_{\lambda\in(\lambda_{lb},\lambda_{ub})}G_Q(\lambda)$\;
		$c^*_{ub}\leftarrow\max(c^*_{ub},G_Q(\lambda^*_{Q,\max}))$\;
		$c^*_{lb}\leftarrow\min(c^*_{lb},G_Q(\lambda^*_{Q,\min}))$\;
	}  
    }
    \caption{Newsvendor under incomplete information: poisson demand. \label{newsvendor_poisson_alg}}
  \end{algorithm}
\end{figure}
The code initially computes Garwood's interval ($\lambda_{lb},\lambda_{ub}$) by exploiting the
relationship between the Poisson distribution and the gamma distribution \citep{Swift09}. Then
it computes the critical fractile $\beta$ and the upper and lower bound for the set $\mathcal{Q}$
of candidate order quantities. Finally, it iterates through the elements of this set to 
compute the upper ($c^*_{ub}$) and lower bound ($c^*_{lb}$) for the estimated cost
associated with these candidate order quantities.

\subsection{Example}\label{sec:example_2}

We consider a simple example involving the Newsvendor problem under Poisson demand. 
In our problem, $h=1$, $p=3$, and the demand follows a $\text{Poisson}(\lambda)$ distribution, 
in which parameter $\lambda$ is unknown. We are
given $10$ samples for the demand, which we may use to determine the optimal order
quantity $Q^*$; these are
$\{51,54,50,45,52,39,52,54,50,40\}$.
The real value for parameter $\lambda$, which is used to generate the samples,
is 50. Accordingly, the optimal order 
quantity $Q^*$ is equal to 55 and provides a cost equal to 9.1222. 

We consider $\alpha=0.9$. 
By using Algorithm \ref{newsvendor_poisson_alg} we compute the 
set of candidate order quantities $\mathcal{Q}=\{50,51,\ldots,57\}$ and
the confidence interval for the estimated cost $(c^*_{lb},c^*_{ub})=(8.6803,14.6220).$
Let us consider a strategy
strategy based on the maximum likelihood estimator. In the case of
the Poisson distribution, this estimator takes the following convenient
form, $\frac{1}{M}\sum_{i=1}^{M}d_i$,
where $d_i$, for $i=1,\ldots,M$ are the observed samples. Therefore,
according to the above samples, the maximum likelihood estimator
for $\lambda$ is 48.7. By using a demand that follows a Poisson
distribution with mean rate $\lambda=48.7$ in Eq. \ref{eq:critical_fractile_disc} we obtain
a candidate optimal order quantity $\widehat{Q}^*=53$ and
an estimated expected cost of $9.0035$. However, such a strategy
does not provide any information on the reliability of the above estimates.
In fact, the actual cost associated with 
this order quantity, when $\lambda=50$, is 9.3693.
Conversely, our approach reports the $\alpha$ confidence
interval $(8.9463,11.0800)$ for the expected cost associated with $\widehat{Q}^*=53$,
which in this case includes the actual cost a decision maker will face in case he decides to order
53 units.

Similar issues occur for the Bayesian approach presented
by \cite{RePEc:eee:ejores:v:98:y:1997:i:3:p:555-562}. For the 
sample set presented above, this approach
suggests ordering 54 units and estimates a cost of 9.4764, 
but it does not provide any 
information on the reliability of these estimates. In contrast,
for an order quantity of 54 units our approach reports the $\alpha$ confidence
interval $(9.0334,10.3374)$ for the expected cost, which include the actual
cost 9.1530 associated with this quantity when $\lambda=50$.

\section{Exponential demand}\label{sec:exp}

A random demand $\text{exp}(\lambda)$ is said to be distributed according to an
exponential law with rate parameter $\lambda>0$ if its probability density function is
\[f(\lambda,k)=\lambda e^{-\lambda k},~~~~~k\geq0;\]
the expected value of $\text{exp}(\lambda)$ is $1/\lambda$.

In the context of the Newsvendor, the exponential distribution may occur in two cases.
An exponentially distributed random variable $\text{exp}(\lambda)$ with rate parameter $\lambda$
can represent the inter-arrival time between two unit demand occurrences in a Poisson
process with rate parameter $\lambda$. Alternatively,
an exponentially distributed random variable $\text{exp}(\lambda)$ can
represent the total demand over the Newsvendor planning horizon. It is clear that the first case can be
easily reduced to the case of a random demand that follows a Poisson distribution
with rate parameter $\lambda$. Such a situation can be handled by following 
the discussion in the previous section.
In the second case, by using Eq. \ref{eq:critical_fractile},  we easily obtain the optimal order quantity 
$Q^*$ for $\text{exp}(\lambda)$. This is simply 
\begin{equation}\label{eq:opt_q_exp}
Q^*=-\frac{1}{\lambda}\ln\left(\frac{h}{h+p}\right).
\end{equation}
We shall give an example. Consider a random demand $\text{exp}(1/50)$ with mean $50$. 
Let $h=1$ and $p=3$, therefore $h/(h+p)=0.25$. The optimal order quantity is $Q^*=69.32$. 
Furthermore, consider the cost function
\[G(Q)=h\int_{0}^{Q}(Q-i)f(\lambda,i)\mathrm{d} i+p\int_{Q}^{\infty}(i-Q)f(\lambda,i)\mathrm{d} i,\]
where $f(\lambda,\cdot)$ denotes the probability density function of $\text{exp}(\lambda)$. Rewrite
\[G(Q)=h(Q-\frac{1}{\lambda})+(h+p)\int_{Q}^{\infty}(1-F(\lambda,i))\mathrm{d} i,\]
where $F(\lambda,\cdot)$ denotes the cumulative distribution function of $\text{exp}(\lambda)$. By noting that 
\begin{equation}\label{eq:opt_cost_exp}
G(Q)=\frac{h+p}{\lambda}\left( \frac{h}{h+p}(\lambda Q-1)+e^{-\lambda Q}\right),
\end{equation}
the optimal cost is $G(Q^*)=69.32$.

We shall now consider, also in this case, the situation in which the parameter $\lambda$ is not known and
the decision maker is given a set of $M$ past realizations of the demand. As in the previous case,
from these realizations he has to infer the range of candidate optimal order quantities and upper and lower bounds for the cost associated with these quantities. This range will contain the real optimum according to a prescribed confidence probability.

\subsection{Confidence intervals for the exponential distribution}

We discuss the exact confidence interval that can be used
to estimate the rate parameter $\lambda$ of the exponential distribution. Consider
a set of $M$ samples $d_i$ drawn from a random variable that is
distributed according to an exponential law with unknown parameter $\lambda$. 
We rewrite $\bar{d}=\sum_{i=0}^M d_i$. Since the sum of $M$ independent and identically distributed exponential 
random variables with rate parameter $\lambda$ is a random variable $\text{gamma}(M,1/\lambda)$ that follows a 
gamma distribution with shape parameter $M$ and scale parameter
$1/\lambda$, the $\alpha$ confidence interval for the unknown parameter 
$\lambda$ is $(\lambda_{lb},\lambda_{ub})$, where
\[
\begin{array}{l}
\lambda_{lb}=\min\{\lambda| \Pr\{\text{gamma}(M,1/\lambda)\geq \bar{d}\}\geq (1-\alpha)/2\},\\
\lambda_{ub}=\max\{\lambda| \Pr\{\text{gamma}(M,1/\lambda)\leq \bar{d}\}\geq (1-\alpha)/2\}.
\end{array}
\]
A closed form expression for this confidence interval 
--- that employs quantiles from the $\chi^2$ distribution ---
was proposed by \citet[chap. 10]{trivedi01} and takes the following form.
Let $\chi^2_{n}$ denote the chi-square distribution 
with $n$ degrees of freedom, and $G^{-1}(\chi^2_{n},\cdot)$ denote the 
inverse cumulative distribution function of $\chi^2_{n}$. We can write
\[\lambda_{lb}=\frac{G^{-1}(\chi^2_{2 M},(1-\alpha)/2)}{2\bar{d}},~~~~~~~~
\lambda_{ub}=\frac{G^{-1}(\chi^2_{2 M},(1+\alpha)/2)}{2\bar{d}}.\]
Furthermore, it is possible to express this interval using quantiles 
from the gamma distribution. More specifically, the lower endpoint 
is the $(1-\alpha)/2$-quantile of a gamma distribution with shape parameter 
$M$ and scale parameter $1/\bar{d}$, and the upper endpoint is the $(1+\alpha)/2$-quantile 
of a gamma distribution with shape parameter $M$ and scale parameter $1/\bar{d}$.

\subsection{Solution method employing statistical estimation based on classical theory of probability}\label{sec:solution_method_exp}

Consider the confidence interval
$\left(\lambda_{lb},\lambda_{ub}\right)$ for the unknown parameter
$\lambda$ of the exponential demand. Let $Q^*_{lb}$ be the optimal
order quantity for the Newsvendor problem under an $\text{exp}(\lambda_{ub})$
demand and $Q^*_{ub}$ be the optimal order quantity for the Newsvendor
problem under an $\text{exp}(\lambda_{lb})$ demand. Recall that
$\lambda$ is a rate, this is the reason why the optimal order quantity
for the Newsvendor problem under an $\text{exp}(\lambda_{lb})$ gives an {\em
upper bound} $Q^*_{ub}$ for the real optimal order quantity. 
Clearly, the optimal order
quantity $Q^*$ lies in the interval $(Q^*_{lb},Q^*_{ub})$.

Let us write the expected total cost associated with an order quantity $Q$
for a given demand rate $\lambda>0$,
\[G_\lambda(Q)=h(Q-\frac{1}{\lambda})+(h+p)\int_{Q}^{\infty}(1-F(\lambda,i))\mathrm{d} i,\]
it is known that this function is convex. Then, consider the function 
\begin{equation}\label{eq:G(lambda)_exp}
G_Q(\lambda)=h(Q-\frac{1}{\lambda})+(h+p)\int_{Q}^{\infty}(1-F(\lambda,i))\mathrm{d} i,
\end{equation}
in which the order quantity $Q$ is fixed and in which we vary the
demand rate $\lambda\geq 0$.  Unfortunately, $G_Q(\lambda)$ is not
convex in the continuous parameter $\lambda$. Nevertheless, we 
prove a number of properties for this function.
\begin{thm}\label{th:theorem_g}
$\lim_{\lambda\rightarrow 0} G_Q(\lambda)=\infty$,
$\lim_{\lambda\rightarrow \infty} G_Q(\lambda)=hQ^-$, the
function admits a single global minimum $\lambda^*$,
it is strictly increasing for $\lambda>\lambda^*$ and strictly 
decreasing for $\lambda<\lambda^*$.
\end{thm}
The proof is given in Appendix \ref{sec:appendix_III}.

Because of the properties of $G_Q(\lambda)$ introduced in Theorem \ref{th:theorem_g}
we can employ a simple line search procedure in order to find 
the $\lambda^*$ that minimizes or maximizes this function over a given interval.

Since the optimal order quantity $Q^*$ is, with confidence probability
$\alpha$, a value in the interval $(Q^*_{lb},Q^*_{ub})$, it is therefore easy to compute
upper ($c^*_{ub}$) and lower ($c^*_{lb}$) bounds for the cost that a manager 
will face, with confidence probability $\alpha$, whatever order quantity he 
chooses in this interval.
\begin{thm}\label{th:theorem_cost_exp}
The lower bound is
\[c^*_{lb}=G_{Q^*_{lb}}(\lambda_{ub})\]
the upper bound is
\[c^*_{ub}=\max\{G_{Q^*_{lb}}(\lambda_{lb}),G_{Q^*_{ub}}(\lambda_{ub})\}.\]
\end{thm}
The proof is given in Appendix \ref{sec:appendix_III}.

Unlike the previous cases, it is not straightforward to extend the above reasoning to 
the case in which unobserved lost sales occurred and the 
$M$ observed past demand data, $d_1,\ldots,d_M$, only reflect the number of customers that
purchased an item when the inventory was positive. This is due to the fact that 
the distribution of the general sum of exponential random variables is not exponential, rather 
it is Hypoexponential. We therefore leave this discussion as a future research
direction.

\subsection{Algorithm}

The computational procedure for exponential demand is presented in Algorithm
\ref{newsvendor_exponential_alg}. 
\begin{figure}[h!]
\vskip 1pc
  \begin{algorithm}[H]
    
    \SetKwInOut{Input}{input}
    \SetKwInOut{Output}{output}
    \Input{$M$ past demand realizations $d_i$, $i=1,\ldots,M$;\\
    the holding cost: $h$;\\
    the penalty cost: $p$;\\
    the confidence probability: $\alpha$.}
    \Output{the set $\mathcal{Q}$ of candidate order quantities;\\
    the interval $(c^*_{lb},c^*_{ub})$ for the estimated cost.}
    \BlankLine
    \Begin{
    	$a\leftarrow M$\;
	$b\leftarrow \sum_{i=1,\ldots,M}d_i$\;
         $\lambda_{ub}\leftarrow$\texttt{InverseCDF[GammaDistribution($a,1/b$),$(1+\alpha)/2$]}\;
         $\lambda_{lb}\leftarrow$\texttt{InverseCDF[GammaDistribution($a,1/b$),$(1-\alpha)/2$]}\;
         $\beta\leftarrow p/(p+h)$\;
         $Q^*_{ub}\leftarrow$\texttt{InverseCDF[PoissonDistribution($\lambda_{ub}$),$\beta$]}\;
         $Q^*_{lb}\leftarrow$\texttt{InverseCDF[PoissonDistribution($\lambda_{lb}$),$\beta$]}\;
         $\mathcal{Q}\leftarrow \{Q^*_{lb},\ldots,Q^*_{ub}\}$\;
    	$c^*_{ub}\leftarrow \max\{G_{Q^*_{lb}}(\lambda_{lb}),G_{Q^*_{ub}}(\lambda_{ub})\}$\;
	$c^*_{lb}\leftarrow G_{Q^*_{lb}}(\lambda_{ub})$\;
    }
    \caption{Newsvendor under incomplete information: exponential demand. \label{newsvendor_exponential_alg}}
  \end{algorithm}
\end{figure}
The code initially computes the confidence interval ($\lambda_{lb},\lambda_{ub}$) by exploiting the
relationship between the exponential distribution and the gamma distribution. Then
it computes the critical fractile $\beta$ and the upper and lower bound for the set $\mathcal{Q}$
of candidate order quantities. Finally, it computes the upper ($c^*_{ub}$) and lower bound ($c^*_{lb}$) f
or the estimated cost associated with these candidate order quantities by exploiting 
Theorem \ref{th:theorem_cost_exp}.
\subsection{Example}\label{sec:example_3}

We consider a simple example involving the Newsvendor problem under exponential demand. 
In our problem, $h=1$, $p=3$, and the demand is a random variable $\text{exp}(\lambda)$ 
for which parameter $\lambda$ is unknown. We are
given $10$ samples for the demand, which we may use to determine the optimal order
quantity $Q^*$; these are
\[\{39.79,39.26,32.21,0.51,107.03,72.87,45.23,20.12,26.46,56.80\}.\]
The real value for parameter $\lambda$, which is used to generate the samples,
is $1/50$. Accordingly, the optimal order 
quantity $Q^*$ is equal to 69.31 and provides a cost equal to 69.31. 

We consider $\alpha=0.9$. The $\alpha$ confidence interval for the 
demand rate $\lambda$ is $(0.0123211,0.0356664)$.
By using Algorithm \ref{newsvendor_exponential_alg} we compute the 
range of candidate order quantities $\mathcal{Q}=(38.86, 112.51)$ and
the confidence interval for the estimated cost $(c^*_{lb},c^*_{ub})=(38.86, 158.81).$
A plot of the expected total cost as a function of $\lambda\in(0.0123211,0.0356664)$ 
and of $Q\in(38.86, 112.51)$ is shown in Appendix \ref{sec:appendix_IV}. 
Let us consider a 
strategy based on the maximum likelihood estimator for the demand rate $\lambda$. In the case of
the exponential distribution, this estimator takes the following
form, $\widehat{\lambda}=\frac{N}{\sum_{i=1}^{N}d_i}$,
where $d_i$, for $i=1,\ldots,N$ are the observed samples. Therefore,
according to the above samples, the maximum likelihood estimator
for $\lambda$ is $\widehat{\lambda}=0.0227099$. By using a
rate $\lambda=0.0227099$ in Eq. \ref{eq:critical_fractile_disc} we obtain
a candidate optimal order quantity $\widehat{Q}^*=61.04$ and
an estimated expected cost of $61.04$.  
As previously remarked, such a strategy
does not provide any information on the reliability of the above estimates.
In fact, the actual cost associated with 
this order quantity, when $\lambda=1/50$, is 70.03.
Conversely, our approach reports the $\alpha$ confidence
interval $(45.71, 132.90)$ for the expected cost associated with $\widehat{Q}^*=61.04$,
which in this case includes the actual cost a decision maker will face in case he decides to order
61.04 units.

Similarly, the Bayesian approach presented
by \cite{RePEc:eee:ejores:v:98:y:1997:i:3:p:555-562}
suggests an order quantity of 59.14 units and estimates a cost of 65.05.
As discussed, this strategy does not provide any 
information on the reliability of these estimates. Conversely,
for an order quantity of 59.14 units, our approach reports the $\alpha$ confidence
interval $(44.71, 134.63)$ for the expected cost, which includes the actual
cost 70.42 associated with this quantity when $\lambda=1/50$.

\section{Discussion and future works}\label{sec:dis}

In this section we first discuss advantages of our strategy, based on Neyman's method of confidence intervals,
with respect to existing frequentist and Bayesian approaches to the Newsvendor problem under
sampled demand information. Secondly, we discuss limitations of
our work and possible future research directions.

\subsection{Comparison with frequentist and Bayesian approaches}
 
Bayesian approaches such as the one proposed by \cite{RePEc:eee:ejores:v:98:y:1997:i:3:p:555-562} 
present several theoretical and practical drawbacks that we are now going
to enumerate. Theoretical drawbacks of the Bayesian approach to parameter estimation
are illustrated by \citet[p. 343]{Ney37}. The first issue raised by Neyman is the fact
that the unknown parameter $a$ of the demand distribution is
not a random variable, therefore assigning a prior or posterior distribution to it
has no meaning. In fact, one may try to employ the prior probability 
distribution --- for instance a uniform distribution in 
Hill's case --- to compute $\Pr\{a < \Delta\}$, where $\Delta\in(0,1)$. 
Of course, $a$ is not a random variable therefore
this probability should be either $1$ or $0$ depending on $\Delta$. Therefore, talking
about prior or posterior distribution for $a$ can only represent an approximation.
It is often stated that in Bayesian probability prior and posterior
distributions are meant to represent a ``state of knowledge'',
that is decision maker's uncertainty about the unknown quantity $a$. However, 
problems arise immediately as soon as one tries to interpret the meaning
of the prior and posterior distribution in light of classical probability theory. 
Neyman, in fact, also points out that, even if the unknown parameter $a$ is
a random variable, the posterior distribution $f(d|b)$ for the random demand $d$, 
computed as illustrated in \cite{RePEc:eee:ejores:v:98:y:1997:i:3:p:555-562}, 
does not generally have the property
serving as a definition of the elementary probability law of the observed data $b$. 
In particular, this distribution is not compatible with the classical definition of probability, in the
sense that, if we repeat an experiment an infinite number of times, 
the observed frequency does not converge to the probability
predicted by such a distribution.
For instance, consider once more the example
presented in Section \ref{sec:example_1}. Instead of 
having a single set of demand observations, we now consider
$M$ experiments, with $M$ large, in each of which we observe 10 demand realizations.
Let $b_i$ be the demand set observed in experiment $i=1,\ldots,M$. 
For each experiment $i$,  we construct the posterior distribution, $f_i(d_i|b_i)$,
of the random demand $d_i$ from $10$ demand observations, according to the
Bayesian strategy discussed by \cite{RePEc:eee:ejores:v:98:y:1997:i:3:p:555-562}. 
Neyman points out that,
in general,  when we select two values $\Delta_1$ and $\Delta_2$,  
and we compute $p_i=\Pr\{\Delta_1<d_i<\Delta_2\}$, 
the quantity $1/M\sum_{i=1}^Mp_i$ will not converge --- as it should, according
to the law of large numbers --- to its real value 
$\Pr\{\Delta_1<\text{bin}(50;q)<\Delta_2\}$, where $q$ is the real value of $a$. 
In our example, we set $\Delta_1=26$ and $\Delta_2=28$. 
$\Pr\{\Delta_1<\text{bin}(50;0.5)<\Delta_2\}=0.1747$.
Nevertheless, if we estimate $p_i=\Pr\{\Delta_1<d_i<\Delta_2\}$ 
in each experiment $i$ by using the posterior 
distribution $f(d_i|b)$ discussed by 
\cite{RePEc:eee:ejores:v:98:y:1997:i:3:p:555-562} for a binomial demand, 
the estimated probability  $1/M\sum_{i=1}^Mp_i$ converges to 0.2654, when the experiment is repeated 
$M=100000$ times. This essentially differs from the real value, as Neyman remarks.
For this reason, using the posterior distribution in place of the original
distribution in the problem of interest is a strategy that may lead to misleading results.

In practice, a direct consequence is that, although for large samples asymptotic results may be obtained (see e.g. \cite{DBLP:journals/rda/BensoussanCRS09}), it becomes hard to 
assess the quality of estimates produced for small sample set. 
It is fairly simple to observe this latter fact by considering once more the example
presented in Section \ref{sec:example_1}. By using the 
approach in \cite{RePEc:eee:ejores:v:98:y:1997:i:3:p:555-562},
we obtain the posterior distribution for the random demand out of the 10 samples, and
then we use this distribution in order to compute an estimate of $Q^*$, which in the
particular example we consider is 29. We also employ the posterior distribution in the
Newsvendor cost function in order to estimate the cost associated with the optimal order quantity
selected; this turns out to be 4.6692. Clearly, Hill's approach does not give any hint on the ``quality'' of the estimates produced, which in this particular case are relatively poor. In particular, based on the available data, we do not know with what frequency 
the order quantity may substantially over- or underestimate the real
optimum order quantity, and how far the prescribed order quantity is likely to be from the 
real optimum order quantity. The same, of course, holds for the estimated optimum cost. 
Our approach based on Neyman's framework, in contrast, suggests that, with 90\% confidence, 
the optimal order quantity --- that is 27 --- lies between 27 and 31,
and that the optimum cost --- that is 4.4946 ---
lies in the interval (4.4268,7.2205). Then, when a given heuristic 
suggests ordering 29 units --- for instance according to Hill's Bayesian approach ---
our approach can be used to derive a 90\% confidence interval for 
the cost associated with this decision, that is  (4.4487, 4.9528). This 
interval actually covers, in this specific case, the real cost associated with the decision of 
ordering 29 units, i.e. 4.8904. In general, the interval will cover the real cost according
to the prescribed confidence probability. Conversely, Hill's approach suggests that the cost 
associated with ordering 29 units is 4.6692, but it provides no indication on the reliability of this estimate. This shows a practical exemplification of how our approach can be employed to effectively complement existing Bayesian approaches under small sample sets. In fact, we must also 
underscore the fact that Bayesian approaches such as the one in 
\citep{RePEc:eee:ejores:v:98:y:1997:i:3:p:555-562} represent
very effective and practical heuristics for order quantity selection. 

Similar issues arise in classical frequentist approaches. For the example in 
Section \ref{sec:example_1}, an approach based on the maximum likelihood 
estimator, as remarked, suggests an order quantity of 29 units. The estimated cost 
according to this strategy is 4.4614. Nevertheless, 
also in this case we have no indication on the reliability of this estimate.
\cite{Kevork_10} derives maximum likelihood estimators for the optimal order 
quantity and for the maximum expected profit.
 The asymptotic distribution for these estimators are then derived and asymptotic 
confidence intervals are extracted for the corresponding true quantities. Unfortunately,
as the author remarks, these intervals are only asymptotically exact and do 
not provide the prescribed confidence, i.e. they are biased, for small sample sets.

By using our approach based on Neyman's framework, the confidence interval
produced for the unknown parameter of the binomial demand is always guaranteed to cover its 
actual value according to the prescribed confidence probability. This probability is controlled by the 
decision maker and influences the size of the interval. Neyman's method 
differs advantageously from the Bayesian approach by being independent 
of a priori information about the unknown parameter. This approach 
remains valid even if the unknown parameter is a random variable. By using our approach,
it is possible to immediately translate the confidence interval
for the unknown parameter into a confidence interval for the
order quantity and for the actual 
cost. Intuitively, in the example presented in Section \ref{sec:example_1}, 
if we observe a set of 10 samples over and over again and we repeat 
our analysis, the intervals produced will cover the real 
optimum order quantity and the associated
cost according to the prescribed probability. In contrast to other existing frequentist of Bayesian approaches our approach provides explicit and exact likelihood guarantees that can be easily interpreted in the context of classical probability theory. Furthermore, by using confidence intervals, the decision maker has
a better control on the risk of exceeding a certain cost and a better outlook on
the range of order quantities that may be optimal according to the observed demands,
especially when a limited set of samples is employed. 

\subsection{Limitations and future works}

Our analysis is limited to three maximum entropy probability distributions 
in the exponential family \citep{citeulike:1512305}, each of which features a single parameter that
must be estimated. As shown by \cite{Harremoes2001},
the binomial and the Poisson are maximum entropy probability distributions for the 
case in which all we know about the distribution of a random demand is that it 
has positive mean and discrete support that goes from 0 to a maximum value $N$ (binomial) or 
to infinity (Poisson). The exponential distribution is the maximum entropy 
probability distribution for the case in which all we know about the 
distribution of a random demand is that it has positive mean and 
continuous support that goes from 0 to infinity. These considerations show
how broadly applicable the results in this work are.
In this work, the normal distribution --- which is part of the exponential family 
and which is also a maximum entropy probability distribution --- has 
not been considered. The analysis on the normal distribution is complicated by the
fact that two parameters, mean and variance, must be considered. Then a number
of cases naturally arise: unknown mean and known variance, unknown variance and
known mean, etc. For this reason, in order to keep the size and the scope of the
discussion limited, we decided to leave this discussion as a future work. 
Furthermore, in principle it may be possible to extend the analysis to other distributions such as the multinomial, for which confidence intervals are surveyed in \citep{citeulike:12601615, citeulike:12601619}; or the Johnson translation system \citep{Johnson1949}, if exact or approximate expressions for the confidence regions of its unknown parameters were available. Unfortunately, we are not aware of any work that investigated these confidence regions.

\section{Conclusions}\label{sec:con}

We considered the problem of controlling the inventory of a single item with stochastic demand over a single period. We introduced a novel strategy to address the issue of demand estimation in single-period inventory optimization problems. Our strategy is based on the theory of statistical estimation. We employed confidence interval analysis in order to identify a range of candidate order quantities that, with prescribed confidence probability, includes the real optimal order quantity for the underlying stochastic demand process with unknown parameter(s). In addition, for each candidate order quantity that is identified, our approach can produce an upper and a lower bound for the associated cost. We applied our novel approach to three demand distribution in the exponential family: binomial, Poisson, and exponential. For two of these distributions we also discussed the
case in which the decision maker faces unobserved lost sales.
Numerical examples are presented in which we showed how our approach complements existing strategies based on maximum likelihood estimators or on Bayesian analysis. 
In particular, we showed that our approach does not provide a 
single order quantity recommendation and a point estimate for the associated cost, 
but --- according to a prescribed confidence level --- a set of candidate 
optimal order quantities and, for each of these, a confidence interval for the associated cost. 
This advanced information can be employed, together with existing frequentist or Bayesian approaches, to better assess the impact of a given decision.

\vspace{-1em}

\section*{Appendix I: proofs of statements for Binomial demand}\label{sec:appendix_I}

Consider Eq. \ref{eq:G(q)} in Section \ref{sec:solution_method}, 
it can be proved that $G_Q(q)$ is
convex in the continuous parameter $q$. Firstly, we rewrite Eq. \ref{eq:G(q)} as
\begin{equation}\label{eq:G(q)_1}
G_Q(q)=h(Q-Nq)+(p+h)\sum_{i=Q}^{N}(1-\Pr\{\text{bin}(N;q)\leq i\}).
\end{equation}
We now show that the second derivative of this function is
positive. Of course, this is equivalent to proving that
\[
\frac{d^2}{dq^2}(p+h)\sum_{i=Q}^{N}(1-\Pr\{\text{bin}(N;q)\leq i\})\geq 0.
\] 
\begin{thm}\label{th:convexity_gq_binomial}
For $Q\leq N$, \[\frac{d^2}{dq^2}\sum_{i=Q}^{N}(1-\Pr\{\text{bin}(N;q)\leq
i\})\] is a positive function of $q\in (0,1)$. 
\end{thm}
\begin{pf}[Theorem \ref{th:convexity_gq_binomial}]
We introduce the following notation:
\[
\begin{array}{l}
f(i;N,q) = \Pr\{\mbox{\it Bin\/}(N;q)=i\} = \binom{N}{i} q^i (1-q)^{N-i},\\
F(i;N,q) = \Pr\{\mbox{\it Bin\/}(N;q) \le i\}.
\end{array}
\]
We reduce the convexity of $G_Q(q)$ to showing that
\[
\frac{d^2}{dq^2}\sum_{i=Q}^N F(i;N,q) \le 0.
\]
Using the regularized incomplete beta function:
\[
F(i;N,q) = (N-i) \binom{N}{i} \int_0^{1-q} t^{N-i-1} (1-t)^i dt;
\]
differentiating under the integral sign by Leibniz's rule:
\[
\begin{array}{rcl}
\frac{d}{dq} F(i;N,q) &=& -(N-i) \binom{N}{i} (1-q)^{N-i-1} q^i\\
&=& -N\binom{N-1}{i}q^i(1-q)^{N-i-1}\\
&=& -N f(i;N-1,q)\\
&=& -N[F(i;N-1,q)-F(i-1;N-1,q)];
\end{array}
\]
using this recursive relationship:
\[
\begin{array}{rcl}
\frac{d^2}{dq^2} F(i;N,q) &=& -N[-(N-1)(F(i;N-2,q)-F(i-1;N-2,q))+\\
&& (N-1)(F(i-1;N-2,q)-F(i-2;N-2,q))]\\
&=& N(N-1)[f(i;N-2,q)-f(i-1;N-2,q)];
\end{array}
\]
and summing over $i$, all terms cancel out except the first and last:
\[
\frac{d^2}{dq^2} \sum_{i=Q}^N F(i;N,q) = N(N-1)[f(N;N-2,q)-f(Q-1;N-2,q)].
\]
However, $f(N;N-2,q)=0$ because it represents the probability of $N$
successes in $N-2$ trials, so
\[
\frac{d^2}{dq^2} \sum_{i=Q}^N F(i;N,q) = -N(N-1)f(Q-1;N-2,q) \le 0
\]
\qed \end{pf}

\section*{Appendix II: proofs of statements for Poisson demand}\label{sec:appendix_II}

Consider Eq. \ref{eq:G(lambda)} in Section \ref{sec:solution_method_poisson}, 
it can be proved that $G_Q(\lambda)$ is convex in the continuous parameter $\lambda$. 
Firstly, we rewrite
Eq. \ref{eq:G(lambda)} as
\begin{equation}\label{eq:G(lambda)_1}
G_Q(\lambda)=h(Q-\lambda)+(h+p)\sum_{i=Q}^{\infty}(1-\Pr\{\text{Poisson}(\lambda)\leq i\}).
\end{equation}
We now show that the second derivative of this function is
positive. Of course, this is equivalent to proving that
\[
\frac{d^2}{d\lambda^2}(h+p)\sum_{i=Q}^{\infty}(1-\Pr\{\text{Poisson}(\lambda)\leq
i\})\geq 0.
\]
Therefore, we have to prove that
\[
\frac{d^2}{d\lambda^2}-\sum_{i=Q}^{\infty}\Pr\{\text{Poisson}(\lambda)\leq i\}\geq 0.
\]
\begin{thm}\label{th:convexity_gq_poisson}
For $Q\geq 0$,
\[
\frac{d^2}{d\lambda^2}-\sum_{i=Q}^{\infty}\Pr\{\text{Poisson}(\lambda)\leq i\}
\] 
is a positive function of $\lambda\geq 0$. 
\end{thm}
\begin{pf}[Theorem \ref{th:convexity_gq_poisson}]
The following derivations prove convexity for the above expression.
\[\frac{d^2}{d\lambda^2}-\sum_{i=Q}^{\infty}\Pr\{\text{Poisson}(\lambda)\leq i\}=\] 
\[\frac{d^2}{d\lambda^2}-\sum_{i=Q}^{\infty}e^{-\lambda}\sum_{k=0}^{i}\frac{\lambda^k}{k!}=\]  
\[-\left(\frac{d^2}{d\lambda^2}e^{-\lambda}\sum_{i=Q}^{\infty}\sum_{k=0}^{i}\frac{\lambda^k}{k!}\right)=\] 
\[-\left(\frac{d}{d\lambda}-e^{-\lambda}\sum_{i=Q}^{\infty}\sum_{k=0}^{i}\frac{\lambda^k}{k!}+
           \frac{d}{d\lambda}e^{-\lambda}\sum_{i=Q}^{\infty}\sum_{k=1}^{i}\frac{\lambda^{k-1}}{(k-1)!}\right)=\] 
           
\[-\left(e^{-\lambda}\sum_{i=Q}^{\infty}\sum_{k=0}^{i}\frac{\lambda^k}{k!}
           -e^{-\lambda}\sum_{i=Q}^{\infty}\sum_{k=1}^{i}\frac{\lambda^{k-1}}{(k-1)!}
           -e^{-\lambda}\sum_{i=Q}^{\infty}\sum_{k=1}^{i}\frac{\lambda^{k-1}}{(k-1)!}+
           e^{-\lambda}\sum_{i=Q}^{\infty}\sum_{k=2}^{i}\frac{\lambda^{k-2}}{(k-2)!}\right)=\] 
           
\[e^{-\lambda}\left(-\sum_{i=Q}^{\infty}\sum_{k=0}^{i}\frac{\lambda^k}{k!}
           +\sum_{i=Q}^{\infty}\sum_{k=1}^{i}\frac{\lambda^{k-1}}{(k-1)!}
           +\sum_{i=Q}^{\infty}\sum_{k=1}^{i}\frac{\lambda^{k-1}}{(k-1)!}
           -\sum_{i=Q}^{\infty}\sum_{k=2}^{i}\frac{\lambda^{k-2}}{(k-2)!}\right)=\] 
           
\[\sum_{i=Q}^{\infty}\left(-e^{-\lambda}\sum_{k=0}^{i}\frac{\lambda^k}{k!}
           +2e^{-\lambda}\sum_{k=1}^{i}\frac{\lambda^{k-1}}{(k-1)!}
           -e^{-\lambda}\sum_{k=2}^{i}\frac{\lambda^{k-2}}{(k-2)!}\right).\]  
           
For convenience, we rewrite this expression as           
           
\[\sum_{i=Q}^{\infty}\left(-\mathrm{CDF}(\text{Poisson}(\lambda),i)
           +2\mathrm{CDF}(\text{Poisson}(\lambda),i-1)
           -\mathrm{CDF}(\text{Poisson}(\lambda),i-2)\right)\]  

\noindent where $\mathrm{CDF}$ denotes the cumulative distribution function. By expanding, we obtain

\[\begin{array}{l}
-\mathrm{CDF}(\text{Poisson}(\lambda),Q)+2\mathrm{CDF}(\text{Poisson}(\lambda),Q-1)-\mathrm{CDF}(\text{Poisson}(\lambda),Q-2)+\\
-\mathrm{CDF}(\text{Poisson}(\lambda),Q+1)+2\mathrm{CDF}(\text{Poisson}(\lambda),Q)-\mathrm{CDF}(\text{Poisson}(\lambda),Q-1)+\\
-\mathrm{CDF}(\text{Poisson}(\lambda),Q+2)+2\mathrm{CDF}(\text{Poisson}(\lambda),Q+1)-\mathrm{CDF}(\text{Poisson}(\lambda),Q)+\ldots=
\end{array}\]

\[\mathrm{CDF}(\text{Poisson}(\lambda),Q-1)-\mathrm{CDF}(\text{Poisson}(\lambda),Q-2)=\]

\[e^{-\lambda}\sum_{k=0}^{Q-1}\frac{\lambda^k}{k!}-e^{-\lambda}\sum_{k=0}^{Q-2}\frac{\lambda^k}{k!}=\] 
\[e^{-\lambda}\left(\sum_{k=0}^{Q-1}\frac{\lambda^k}{k!}-\sum_{k=0}^{Q-2}\frac{\lambda^k}{k!}\right)=\] 
\[\frac{e^{-\lambda}\lambda^{Q-1}}{(Q-1)!}\geq0\]
\qed \end{pf}

\section*{Appendix III: proofs of statements for exponential demand}\label{sec:appendix_III}

In this section we provide the proofs for the two theorems introduced in Section \ref{sec:solution_method_exp}.

\begin{pf}[Theorem \ref{th:theorem_g}]
The first fact, that is $\lim_{\lambda\rightarrow 0}
G_Q(\lambda)=\infty$ can be easily verified by simple algebraic
derivations.  We shall therefore prove that $\lim_{\lambda\rightarrow
\infty} G_Q(\lambda)=hQ^-$ and that the function admits a single global
minimum.  Let us split Eq. \ref{eq:opt_cost_exp} into two parts
\begin{equation}\label{eq:G(lambda)_1split}
G_Q(\lambda)=\frac{h+p}{\lambda}e^{-\lambda Q}+
\left(hQ-\frac{h}{\lambda}\right).
\end{equation}
We shall consider the first term
\begin{equation}\label{eq:G(lambda)_2}
\frac{h+p}{\lambda}e^{-\lambda Q},
\end{equation}
and the second term
\begin{equation}\label{eq:G(lambda)_3}
hQ-\frac{h}{\lambda},
\end{equation}
on the right hand side of Eq. \ref{eq:opt_cost_exp}, separately. 

Firstly, we observe that, when $\lambda \rightarrow
\infty$, $G_Q(\lambda)\rightarrow hQ$ from below, that is
$\lim_{\lambda\rightarrow \infty} G_Q(\lambda)=hQ^-$. 
This is due to
the fact that (\ref{eq:G(lambda)_2}) approaches zero faster than
$h/\lambda$ does, i.e.
\[
\lim_{\lambda\rightarrow\infty} \frac{\frac{h+p}{\lambda}e^{-\lambda
Q}}{\frac{h}{\lambda}}=0.
\]
From this fact we immediately infer that the derivative of
$G_Q(\lambda)$ must be equal to zero for at least one value $\lambda$
other than infinity. Furthermore, the derivative of
(\ref{eq:G(lambda)_2})
is negative, strictly increasing for $\lambda > 0$.
The derivative of (\ref{eq:G(lambda)_3})
is positive strictly decreasing for $\lambda > 0$. Therefore there
exists only a single value of $\lambda$ for which the derivative of
(\ref{eq:G(lambda)_2}) and the derivative of (\ref{eq:G(lambda)_3})
add up to zero. This immediately implies that $G_Q(\lambda)$ admits a
single global minimum, it is strictly increasing for
$\lambda>\lambda^*$ and strictly decreasing for $\lambda<\lambda^*$
\qed \end{pf}

\begin{pf}[Theorem \ref{th:theorem_cost_exp}]
Firstly, let us consider $c^*_{lb}$. By definition, this is the expected total cost associated with
the optimal order quantity $Q^*_{lb}$ for the largest possible value $\lambda_{ub}$ that the demand rate 
takes in the confidence interval. Consider a demand rate $\lambda$ and the associated
optimal order quantity $Q^*_{\lambda}$. By substituting $Q$ in Eq. \ref{eq:opt_cost_exp}
with the expression of the optimal order quantity in Eq. \ref{eq:opt_q_exp} we immediately see
that the expected total cost associated with
an optimal order quantity $Q^*_{\lambda}$ is decreasing in the respective demand rate $\lambda$ --- i.e. it is
increasing w.r.t. the expected value $1/\lambda$ of the demand --- it immediately
follows that there exists no other pair $\langle Q^*_{\lambda},\lambda\rangle$, where 
$\lambda\in(\lambda_{lb},\lambda_{ub})$ that ensures a lower expected total cost.

Let $\lambda\in(\lambda_{lb},\lambda_{ub})$ and $Q\in(Q_{lb},Q_{ub})$.
Consider a point in the two dimensional space $\lambda\times Q$, for which
$\lambda=\bar{\lambda}$ and $Q=\bar{Q}$. For any of such points, two cases
can be observed, that is (i)  $\bar{Q}>Q^*_{\bar{\lambda}}$, or (ii) 
$\bar{Q}<Q^*_{\bar{\lambda}}$. A strict equality can be reduced to any of these
two cases. If we are in case (i), then $G_{\bar{\lambda}}(\bar{Q})<G_{\bar{\lambda}}(Q_{ub})$
because of Theorem \ref{th:theorem_g}. $\bar{Q}$ was already an order quantity
larger than the optimal one, therefore $Q_{ub}$ is also an order quantity
larger than the optimal one for a demand rate $\bar{\lambda}$. Consequently,
if we increase the demand rate $\lambda$ (i.e. we decrease the expected demand $1/\lambda$) our
cost can only increase; this means that $G_{\bar{\lambda}}(Q_{ub})<G_{\lambda_{ub}}(Q_{ub})$.
If we are in case (ii), then $G_{\bar{\lambda}}(\bar{Q})<G_{\bar{\lambda}}(Q_{lb})$
because of Theorem \ref{th:theorem_g}. $\bar{Q}$ was already an order quantity
smaller than the optimal one, therefore $Q_{lb}$ is also an order quantity
smaller than the optimal one for a demand rate $\bar{\lambda}$. Consequently,
if we decrease the demand rate $\lambda$ (i.e. we increase the expected demand $1/\lambda$) our
cost can only increase; this means that $G_{\bar{\lambda}}(Q_{lb})<G_{\lambda_{lb}}(Q_{lb})$.
Therefore, the maximum cost, when we let $\lambda$ vary in $(\lambda_{lb},\lambda_{ub})$
and $Q$ vary in $(Q^*_{lb},Q^*_{ub})$, can be either observed at $\langle Q^*_{lb},\lambda_{lb} \rangle$
or at $\langle Q^*_{ub},\lambda_{ub} \rangle$
\qed \end{pf}

\section*{Appendix IV: plot for the expected total cost of the example in Section \ref{sec:example_3}}\label{sec:appendix_IV}
\begin{figure}[h!]
\centering 
\includegraphics[type=eps,ext=.eps,read=.eps,width=0.9\columnwidth]{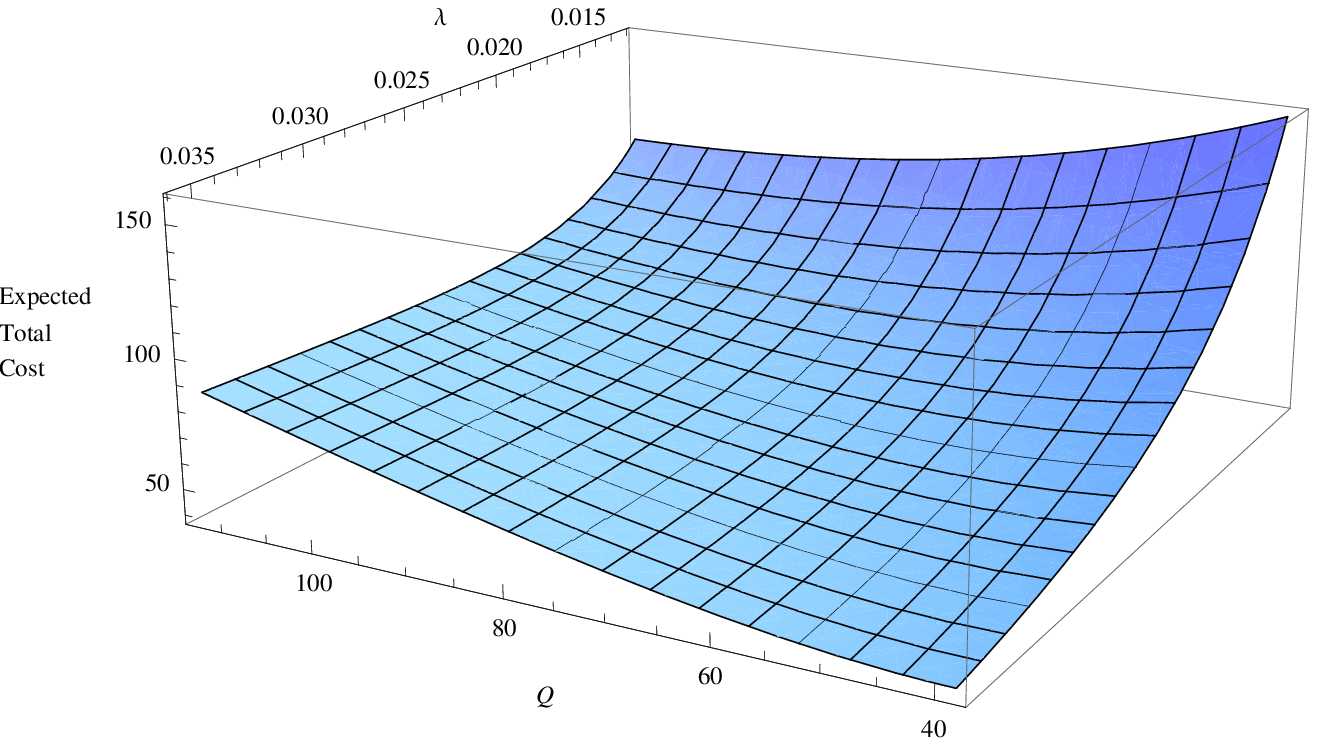}
\caption{Expected total cost as a function of $\lambda\in(\lambda_{lb},\lambda_{ub})\equiv(0.0123211,0.0356664)$
and of $Q\in(Q^*_{lb},Q^*_{ub})\equiv(38.86, 112.51)$. Note that $c^*_{lb}=G_{Q^*_{lb}}(\lambda_{ub})$ and that
$c^*_{ub}=\max\{G_{Q^*_{lb}}(\lambda_{lb}),G_{Q^*_{ub}}(\lambda_{ub})\}=
G_{Q^*_{lb}}(\lambda_{lb})$.
}
\label{fig:newsvendor_exp}
\end{figure}

In Fig. \ref{fig:newsvendor_exp} we provide a graphical outlook of the cost function discussed in Section \ref{sec:example_3}.

\bibliographystyle{plainnat}
\bibliography{newsvendor}

\end{document}